\documentclass[11pt,a4paper,leqno]{article}
\usepackage{a4wide}
\usepackage{latexsym}
\usepackage{amsmath}
\usepackage{amssymb}
\usepackage{stackrel}  
\usepackage{color} 
\usepackage{graphicx}
\usepackage[linesnumbered,ruled,vlined]{algorithm2e}

\usepackage{authblk}

\usepackage{hyperref}

\newtheorem{defin}{Definition}

\newtheorem{prop}{Proposition}
\newtheorem{theo}{Theorem}
\newtheorem{corol}{Corollary}

\newenvironment{proof}{\medskip\par\noindent{\bf Proof}}{\hfill $\Box$
\medskip\par}

\newcommand{\C}{\mathbb{C}}
\newcommand{\N}{\mathbb{N}}
\newcommand{\R}{\mathbb{R}}

\begin{document}
\title{Gevrey versus $q-$Gevrey asymptotic expansions for some linear $q-$difference-differential Cauchy problem}

\author[1]{Alberto Lastra}
\author[2]{St\'ephane Malek}
\affil[1]{Universidad de Alcal\'a, Dpto. F\'isica y Matem\'aticas, Alcal\'a de Henares, Madrid, Spain. {\tt alberto.lastra@uah.es}}
\affil[2]{University of Lille, Laboratoire Paul Painlev\'e, Villeneuve d'Ascq cedex, France. {\tt stephane.malek@univ-lille.fr}}


\date{}

\maketitle
\thispagestyle{empty}
{ \small \begin{center}
{\bf Abstract}
\end{center}

The asymptotic behavior of the analytic solutions of a family of singularly perturbed $q$-difference-differential equations in the complex domain is studied. 

Different asymptotic expansions with respect to the perturbation parameter and to the time variable are provided: one of Gevrey nature, and another of mixed type Gevrey and $q-$Gevrey.

This asymptotic phenomena is observed due to the modification of the norm established on the space of coefficients of the formal solution. The techniques used are based on the adequate path deformation of the difference of two analytic solutions, and the application of several versions of Ramis-Sibuya theorem.

\smallskip

}
\bigskip \bigskip

\section{Introduction}

In the present study, a family of singularly perturbed linear $q$-difference-differential equations of the form
\begin{equation}\label{epralintro}
P(\epsilon^kt^{k+1}\partial_t)\partial_z^{S}u(t,z,\epsilon)=\mathcal{P}(t,z,\epsilon,\partial_t,\partial_z,\sigma_q)u(t,z,\epsilon),
\end{equation}
under initial data 
\begin{equation}\label{cdintro}
(\partial_z^ju)(t,0,\epsilon)=\varphi_j(t,\epsilon),\quad 0\le j\le S-1,
\end{equation}
is studied. In the previous problem, $\epsilon$ acts as a small complex perturbation parameter, and $\sigma_q$ stands for the dilation operator on $t$ variable defined by $\sigma_qf(t)=f(qt)$, for some fixed $q>1$. In (\ref{epralintro}), $S,k$ are positive integers, $P(\tau)\in\C[\tau]$, and the symbol $\mathcal{P}(t,z,\epsilon,\tau_1,\tau_2,\tau_3)$ is a polynomial in $(t,z,\tau_1,\tau_2,\tau_3)$ with holomorphic coefficients on some neighborhood of the origin with respect to the perturbation parameter. The functions $\varphi_j(t,\epsilon)$ are polynomials with respect to $t$ variable, with holomorphic coefficients on some neighborhood of the origin in $\epsilon$. It is worth mentioning the irregular nature of the differential operator in both $P(\epsilon^kt^{k+1}\partial_t)$, and inside $\mathcal{P}$ in the form of polynomial operators in $\epsilon^kt^{k+1}\partial_t$. The precise shape of the main problem and the concrete assumptions considered in it are determined in detail in Section~\ref{sec21}. 

There is an increasing interest on the study of the asymptotic behavior of the solutions to $q-$difference-differential equations in the complex domain. This is the case of the recent works by H. Tahara~\cite{tahara}, H. Yamazawa~\cite{ya1} and H. Tahara and H. Yamazawa~\cite{taya2}; the authors and J. Sanz~\cite{lamasa0} and with T. Dreyfus~\cite{drelasmal}. A different approach via Nevalinna theory is developped in~\cite{hongyanlianzhong,tuyuan}. The importance of applications of $q$-difference equations in the knowledge of wavelets or tsunami and rogue waves is evidenced in recent advances in the field such as~\cite{prrasp,prrasp1}.

In the present work, we construct solutions $u(t,z,\epsilon)$ of the main problem (\ref{epralintro}), (\ref{cdintro})  which are bounded holomorphic functions defined on $\mathcal{T}\times D_R\times \mathcal{E}$, where $\mathcal{T}$ and $\mathcal{E}$ stand for finite sectors of the complex plane with vertex at the  origin, and $D_{R}$ stands for the open disc of radius $R>0$. The main purpose of the present work is to show that the function $\mathcal{E}\ni\epsilon\mapsto u(t,z,\epsilon)$ shows different asymptotic expansions with respect to the perturbation parameter $\epsilon$ when modifying the norm considered in the space of coefficients, the space of holomorphic and bounded functions defined on $\mathcal{T}\times D_R$. The symmetric situation with respect to $\mathcal{T}\ni t\mapsto u(t,z,\epsilon)$ is also considered.

The appearance of such phenomena is due to the existence of small divisors involved in the main equation, in contrast to~\cite{ma17}. In the seminal work~\cite{ma17}, the second author studies the asymptotic solutions with respect to $\epsilon$ of equations of the form
\begin{equation}\label{e98}
Q(\partial_z)u(qt,z,\epsilon)=\sum_{\ell=1}^{D}\epsilon^{\Delta_{\ell}}t^{d_\ell}
c_{\ell}(t,z,\epsilon)R_{\ell}(\partial_z)u(q^{\delta_\ell}t,z,\epsilon)+f(t,z,\epsilon),
\end{equation}
where $q>1$, $D\ge 2$, $\delta_{\ell},\Delta_{\ell},d_{\ell}$ are non-negative integers for $1\le \ell\le D$. The functions $c_{\ell}(t,z,\epsilon)$ for $1\le \ell\le D$ and $f(t,z,\epsilon)$ are bounded holomorphic functions on $D_{r}\times \{z\in\C: |\hbox{Im}(z)|\le \beta'\}\times D_{\epsilon_0}$, for some $r,\beta',\epsilon_0>0$. $Q(\tau)$ and $R_{\ell}(\tau)$ for $1\le \ell\le D$ are polynomials with complex coefficients. All these elements are subjected to further hypotheses not mentioned here for the sake of simplicity.

In~\cite{ma17}, the technique used to solve asymptotically the problem is to construct the analytic solutions to (\ref{e98}) in the form of a $q-$Laplace transform of order $k$, for some $k>0$ which depends on the elements of the problem. This causes the absence of two distinguished asymptotic expansions (even when modifying the norm considered for the function spaces involving the variables $(t,z)$) in contrast to the present situation. In that work, no small divisor appear in the main problem under study, all the asymptotic expansions obtained there being of $q-$Gevrey type with respect to $\epsilon$.

The main inspiration of the present study is~\cite{carrillolastra}, where the authors deal with the formal solutions to systems of dimension $N\ge1$ of the form
\begin{equation}\label{e116}
\epsilon^{\alpha}x^{p}\sigma_{q,x}(y)(x,\epsilon)=F(x,\epsilon,y),
\end{equation}
where $p,\alpha$ are non-negative integers, $q\in\C$ is such that $|q|>1$, $\alpha>0$, and $F(x,\epsilon,y)$ is analytic on some neighborhood of the origin in $\C\times\C\times\C^{N}$, with $F(0,0,0)=0$ and $DF_y(0,0,0)$ being an invertible matrix. Indeed, the unique formal solution of (\ref{e116})
$$\hat{y}(x,\epsilon)=\sum_{n=0}^{\infty}y_n(\epsilon)x^n=\sum_{n=0}^{\infty}u_n(x)\epsilon^n$$
is such that 
\begin{itemize}
\item[(1)] for $p>0$, all $y_n$ converge on some common neighorhood of the origin, whereas $u_n$ converges on the disc of radius $r/|q|^{\left\lfloor n/\alpha\right\rfloor}$ for some $r>0$ and there exist $C=C(q),A=A(q)>0$ such that for all $n\ge0$ one has
$$\sup_{|\epsilon|\le r}|y_n(\epsilon)|\le CA^n|q|^{\frac{n^2}{2p}},\quad \sup_{|\epsilon|\le r/|q|^{\left\lfloor n/\alpha\right\lfloor}}|u_n(x)|\le CA^n.$$
\item[(2)] for $p=0$, $y_n$ and $u_n$ converge on the disc of radius $r/|q|^{ n/\alpha}$ and $r/|q|^{\left\lfloor n/\alpha\right\rfloor}$, respectively, for some $r>0$, and there exist $C=C(q),A=A(q)>0$ such that for all $n\ge0$ one has
$$\sup_{|\epsilon|\le r/|q|^{n/\alpha}}|y_n(\epsilon)|\le CA^n,\quad \sup_{|\epsilon|\le r/|q|^{\left\lfloor n/\alpha\right\lfloor}}|u_n(x)|\le CA^n.$$
\end{itemize}

Observe in the previous result that the coefficients of the formal solutions might be defined in shrinking neighborhoods of the origin, determining power series which have null radius of convergence. 

The procedure followed to solve (\ref{epralintro}), (\ref{cdintro}) analytically is to search for solutions in the form of a Laplace transform of order $k$ (see Section~\ref{secannex}) which transforms the main problem into an auxiliary convolution equation, whose analytic solution satisfies appropriate bounds in order to recover an analytic solution to the main equation via Laplace transform (see Proposition~\ref{prop5}). Sharp bounds satisfied by the solutions to the auxiliary convolution equation are also available, leading to the construction of a finite family of analytic solutions to (\ref{epralintro}), (\ref{cdintro}), say $(u_{p,1}(t,z,\epsilon))_{0\le p\le \varsigma_1-1}$, for some integer $\varsigma_1\ge2$, with $u_{p,1}\in\mathcal{O}_b(\mathcal{T}\times D\times \mathcal{E}_p)$ for all $0\le p\le \varsigma_1-1$. Here, $\mathcal{T}$ stands for a bounded sector in $\C$ with vertex at the origin, $D$ is a neighborhood of the origin, and $\mathcal{E}_p$ is a bounded sector with vertex at the origin belonging to a good covering in $\C^{\star}$ (see Definition~\ref{defin-goodcovering}). Another finite family of analytic solutions to 
(\ref{epralintro}), (\ref{cdintro}), say $(u_{p,2}(t,z,\epsilon))_{0\le p\le \varsigma_2-1}$, for some integer $\varsigma_2\ge2$, is also constructed. For every $0\le p\le\varsigma_2-1$, the solution $u_{p,2}(t,z,\epsilon)$ remains analytic on $\mathcal{T}_p\times D\times \mathcal{E}$, where $\mathcal{E}$ is some bounded sector in $\C$ with vertex at the origin, and so it is $\mathcal{T}_p$, which is an element of a good covering in $\C^{\star}$.

The main results of the present work determine the asymptotic behavior of the two families of  analytic solutions from two radically different topological points of view. It is proved in Theorem~\ref{teo2} the existence of a formal power series in the perturbation parameter, with coefficients in some Banach space of functions which asymptotically approaches each of the analytic solutions in $(u_{p,1}(t,z,\epsilon))_{0\le p\le \varsigma_1-1}$. The asymptotic approximation is measured by means of a $L_1-q-$relative-sup-norm. Such norm is defined on a larger set of formal power series in one of their variables with coefficients being holomorphic functions on some shrinking neighborhood of the origin (see Definition~\ref{defi243}). Under this measurement, the asymptotic behavior is of Gevrey nature (see (\ref{e344})). On the other hand, when incorporating the classical $L_1$-sup norm in the asymptotic approximation, then mixed Gevrey and $q-$Gevrey asymptotic expansions emerge, as it is proved in Theorem~\ref{teo3}. We recall that previous results in the field have also observed such multiscaled asymptotics, such as~\cite{lama17}. Theorem~\ref{teo2} and Theorem~\ref{teo3} are put forward in a symmetric manner regarding time variable, leading to Gevrey and $q$-Gevrey asymptotic relations for the analytic solutions $(u_{p,2}(t,z,\epsilon))_{0\le p\le \varsigma_2-1}$, in Theorem~\ref{teo4} and Theorem~\ref{teo5}. The technique used in the preceeding results leans on the application of the classical version of the so-called Ramis-Sibuya theorem (Theorem~\ref{rst} (RS) in Section~\ref{secRS}) and a $q-$analog of Ramis-Sibuya theorem (Theorem~\ref{trsq} ($q-$RS) in Section~\ref{secRS}).

In brief, the work states different asymptotic expansions with respect to $\epsilon$ and $t$ regarding different sets of analytic solutions to the main problem under study (\ref{epralintro}), (\ref{cdintro}): one of Gevrey order $1/k$ and another of mixed type Gevrey and $q-$Gevrey, when modifying the norm set on the space of coefficients of the formal solution. In addition to this, Gevrey order expansions of order $1/k$ have been observed in both variables $t$ and $\epsilon$ by setting appropriate norms on the spaces of holomorphic functions involved.

The paper is structured as follows. Section~\ref{sec21} is devoted to precise the main problem under study. In the next subsections, we provide different families of analytic solutions (Theorem~\ref{teo1} and Theorem~\ref{teo1bis}) by fixing concise geometries in the problem. The first main results on the asymptotic behavior of the previous families of analytic solutions are stated in Section~\ref{secasI} (Theorem~\ref{teo2} and Theorem~\ref{teo3}) by determining different norms in the space of coefficients of the formal solution. Symmetric results regarding the time variable (Theorem~\ref{teo4} and Theorem~\ref{teo5}) are stated in Section~\ref{secasII}. The work concludes with two annex which complete known facts about Laplace transform and its main properties and several versions of Ramis-Sibuya type theorems, appealed in the paper.

\vspace{0.4cm}

\noindent\textbf{Notation:}

We write $\N:=\{1,2,3,\ldots\}$ and $\N_0:=\N\cup\{0\}$.

For every $r>0$ and $z_0\in\C$, we write $D(z_0,r)$ for the open disc centered at $z_0$ and radius $r$, and for simplicity we denote $D_{r}:=D(0,r)$.

Given a nonempty open set $U\subseteq\C$, and a complex Banach space $\mathbb{E}$, $\mathcal{O}_b(U,\mathbb{E})$ stands for the set of bounded holomorphic functions $h:U\to\mathbb{E}$ which determines a Banach space with the norm of the supremum. For simplicity, we write $\mathcal{O}_b(U)$ instead of $\mathcal{O}_b(U,\C)$. We also denote the set of formal power series in the variable $z$ and coefficients in $\mathbb{E}$ by $\mathbb{E}[[z]]$.

\section{Statement of the main problem and analytic solution}

In this section, we state the main Cauchy problem under study in the present work and provide analytic solutions to it in adequate domains. The procedure heavily rests on the study developed in~\cite{malek22}, and therefore most of the details are omitted, in order to enhance the main purpose of the present study, i.e. to focus on the appearance of different related asymptotic occurrence related to such analytic solutions. 

\subsection{Statement of the main problem}\label{sec21}

Let $k,S$ be positive integers. Fix $q>1$, and $\epsilon_0>0$.

Let $\mathcal{A}\subseteq\N_0^4$ be a finite set, and  $\Delta_{\underline{\ell}}\in\N$ for every $\underline{\ell}\in\mathcal{A}$. We assume that for every $\underline{\ell}=(\ell_0,\ell_1,\ell_2,\ell_3)\in\mathcal{A}$
\begin{equation}\label{e1}
\ell_2<S, \qquad S\ge \ell_2+\ell_3,\qquad \Delta_{\underline{\ell}}\ge \ell_0,
\end{equation}

and fix a polynomial $c_{\underline{\ell}}(z,\epsilon)=\sum_{h\in I_{\underline{\ell}}}c_{\underline{\ell},h}(\epsilon)z^h\in \mathcal{O}_b(D_{\epsilon_0})[z]$, for some finite set $I_{\underline{\ell}}\subseteq\N$. We assume that $c_{\underline{\ell}}(0,\epsilon)\equiv0$. The following hypotheses are also fulfilled: there exists $\Delta\ge 1/2$ such that
\begin{multline}\label{e2}
2(\ell_2-h)\Delta +\ell_0+k\ell_1-2(S-1)\Delta>0,\\ 
\Delta\max\{0,2(\ell_2-h)-1\}-(-h+\ell_2)^2\Delta<\min_{a\in \{S-1,S\}}a(\ell_0+k\ell_1)-a^2\Delta,\\
-2\Delta\ell_3+\ell_0+k\ell_1>0,
\end{multline} 
for all $h\in I_{\underline{\ell}}$ and $\underline{\ell}=(\ell_0,\ell_1,\ell_2,\ell_3)\in\mathcal{A}$. 

\noindent\textbf{Remark:} Less restrictive conditions than those appearing in (\ref{e2}) can be assumed. We have decided to adopt (\ref{e2}) for the sake of simplicity. Observe moreover that the assumptions (\ref{e1}) and (\ref{e2}) are compatible whenever $\ell_0+k\ell_1$ for $(\ell_0,\ell_1,\ell_2,\ell_3)\in\mathcal{A}$ is large enough, together with suitable choices of the remaining constants.

\vspace{0.3cm}

For every $0\le j\le S-1$, we consider a polynomial in time variable, say $\varphi_j(t,\epsilon)\in(\mathcal{O}_b(D_{\epsilon_0}))[t]$, such that $\varphi_j(0,\epsilon)\equiv0$ for $0\le j\le S-1$.

Let $P(\tau)\in\C[\tau]$ with $P(0)\neq0$, and such that there exists $k_1>0$ with
\begin{equation}\label{e3}
k\hbox{deg}(P)\ge k\ell_1+\ell_0+2k_1\ell_3\log(q),\quad\hbox{for } \underline{\ell}=(\ell_0,\ell_1,\ell_2,\ell_3)\in\mathcal{A}.
\end{equation}

The main problem under consideration is the following singularly perturbed linear Cauchy problem
\begin{equation}\label{epral}
P(\epsilon^k t^{k+1}\partial_t)\partial_z^Su(t,z,\epsilon)=\sum_{\underline{\ell}=(\ell_0,\ell_1,\ell_2,\ell_3)\in\mathcal{A}}\epsilon^{\Delta_{\underline{\ell}}}c_{\underline{\ell}}(z,\epsilon)t^{\ell_0}\left((\epsilon^{k}t^{k+1}\partial_t)^{\ell_1}\partial_z^{\ell_2}u\right)(q^{\ell_3}t,z,\epsilon),
\end{equation}
with Cauchy data
\begin{equation}\label{epralCD}
(\partial_z^ju)(t,0,\epsilon)=\varphi_{j}(t,\epsilon),\qquad 0\le j\le S-1.
\end{equation}

\subsection{Construction of analytic solutions to the main problem}

The strategy to find analytic solutions to (\ref{epral}), (\ref{epralCD}) is to search for functions in the form of a Laplace transform of order $k$ along well chosen directions $\gamma\in\R$ to be determined. More precisely, we search for solutions
$$u(t,z,\epsilon)=k\int_{L_{\gamma}}\omega(u,z,\epsilon)\exp\left(-\left(\frac{u}{\epsilon t}\right)^{k}\right)\frac{du}{u},$$
for some suitable function $\omega$ under suitable growth properties at infinity regarding its first variable. 

In view of the properties displayed in Proposition~\ref{propLaplace2} and Proposition~\ref{propLaplace3} on Laplace transform, we reduce the problem to solve the following auxiliary Cauchy problem

\begin{multline}
\partial_z^S\omega(u,z,\epsilon)=\sum_{\underline{\ell}=(\ell_0,\ell_1,\ell_2,\ell_3)\in\mathcal{A},\ell_0=0}\epsilon^{\Delta_{\underline{\ell}}}\frac{c_{\underline{\ell}}(z,\epsilon)}{P(ku^k)}(k(q^{\ell_3}u)^{k_1})^{\ell_1}(\partial_z^{\ell_2}\omega)(q^{\ell_3}u,z,\epsilon)\\
+\sum_{\underline{\ell}=(\ell_0,\ell_1,\ell_2,\ell_3)\in\mathcal{A},\ell_0\ge1}\epsilon^{\Delta_{\underline{\ell}}-\ell_0}\frac{c_{\underline{\ell}}(z,\epsilon)}{P(ku^k)}\frac{u^k}{\Gamma\left(\frac{\ell_0}{k}\right)}\int_0^{u^k}(u^k-s)^{\frac{\ell_0}{k}-1}(k(q^{\ell_3}s^{1/k})^k)^{\ell_1}(\partial_z^{\ell_2}\omega)(q^{\ell_3}s^{1/k},z,\epsilon)\frac{ds}{s}\label{eaux},
\end{multline}
with Cauchy data
\begin{equation}\label{eauxCD}
(\partial_z^j\omega)(u,0,\epsilon)=P_{j}(u,\epsilon),\qquad 0\le j\le S-1,
\end{equation}
where $P_j$ for $0\le j\le S-1$ is determined from $\varphi_j$ from the properties of Laplace transform (see Proposition~\ref{propLaplace1} and Corollary~\ref{coroLaplace1}). Indeed, one can express $\varphi_j$ as the Laplace transform of order $k$ of some polynomial $P_j(u,\epsilon)\in(\mathcal{O}_b(D_{\epsilon_0}))[u]$. More precisely, one has that $P_j(u,\epsilon)=\sum_{h\in J_j}p_{j,h}(\epsilon)u^h$, for some finite subset $J_j\subseteq \N$, and all $0\le j\le S-1$. In this situation, $\varphi_j(t,\epsilon)=\sum_{h\in J_j}\Gamma(h/k)p_{j,h}(\epsilon)t^h$, for all $0\le j\le S-1$, where $\Gamma(\cdot)$ stands for Gamma function.

A finite family of analytic solutions to (\ref{eaux}), (\ref{eauxCD}) is now constructed from the following associated geometric configuration to the problem, on adequate domains.

The following result is a parametric version of Proposition 5 in~\cite{malek22}, whose proof can be adapted without any difficulty to this setting. We only give a sketch of the points to be adapted in our framework.

\begin{prop}\label{prop5}
Let $\varsigma\ge 2$ be an integer. We consider a set $\underline{\mathcal{U}}=(\mathcal{U}_{p})_{0\le p\le \varsigma-1}$ of unbounded sectors with vertex at the origin which satisfies that all the roots of the polynomial $u\mapsto P(ku^k)$ belong to $\C\setminus\left(\bigcup_{p=0}^{\varsigma-1}\mathcal{U}_p\right)$. Assume conditions (\ref{e1}), (\ref{e2}) and (\ref{e3}) hold. Then, there exists $R>0$ such that for every $0\le p\le\varsigma-1$ a series
\begin{equation}\label{eA}
\omega_p(u,z,\epsilon)=\sum_{n\ge0}\omega_{p,n}(u,\epsilon)\frac{z^n}{n!}
\end{equation}
can be constructed determining a holomorphic function on $\mathcal{U}_p\times D_R\times D_{\epsilon_0}$, which solves the auxiliary Cauchy problem (\ref{eaux}), (\ref{eauxCD}). Moreover, there exist $C_3>0$, $u_0>1$, $\alpha\ge 0$ and $k_1>0$ stated in (\ref{e3}) with
\begin{equation}\label{e186}
|\omega_p(u,z,\epsilon)|\le 2C_3|u|\exp\left(k_1\log^2(|u|+u_0)+\alpha\log(|u|+u_0)\right),
\end{equation}
for all $u\in\mathcal{U}_p$, $z\in D_{R}$ and $\epsilon\in D_{\epsilon_0}$. In addition to this, each function $\omega_{p,n}(u,\epsilon)$ satisfies the following properties:
\begin{itemize}
\item[a.] For every $0\le p\le \varsigma-1$ and $n\ge0$, the function $(u,\epsilon)\mapsto \omega_{p,n}(u,\epsilon)$ belongs to $\mathcal{O}(\mathcal{U}_p\times D_{\epsilon_0})$ with
\begin{equation}\label{e194}
\sup_{\epsilon\in D_{\epsilon_0}}|\omega_{p,n}(u,\epsilon)|\le C_3\frac{1}{(2R)^n}n!|u|\exp\left(k_1\log^2(|u|+u_0)+\alpha\log(|u|+u_0)\right),
\end{equation}
for all $u\in\mathcal{U}_p$.
\item[b.] For every $n\ge0$ there exists an analytic function $\omega_n(u,\epsilon)$, which is a common analytic extension of $(u,\epsilon)\mapsto \omega_{p,n}(u,\epsilon)$ for all $0\le p\le \varsigma-1$. This analytic continuation is defined on $D_{R_n}\times D_{\epsilon_0}$ with $R_n=R_0/q^n$, for some small enough $R_0>0$ such that 
$$D_{R_0}\cap\{u\in\C:P(ku^k)=0\}=\emptyset.$$
Moreover, there exist $C_1,C_2>0$ such that
\begin{equation}\label{e195}
\sup_{\epsilon\in D_{\epsilon_0}}|\omega_{n}(u,\epsilon)|\le C_1(C_2)^n\frac{n!}{q^{n^2\Delta}}|u|,
\end{equation}
for $u\in D_{R_n}$, where $\Delta$ is introduced in (\ref{e2}).
\item[c.] For every $0\le p\le \varsigma-1$ and $n\ge0$, the function $(u,\epsilon)\mapsto \omega_{p,n}(u,\epsilon)$ is bounded holomorphic on the sectorial annulus
$$\mathcal{A}_{p,h}=\left\{u\in\mathcal{U}_p:\frac{R_0}{q^{h+1}}\le |u|\le\frac{R_0}{q^{h}}\right\}$$
for all $0\le h\le n-1$, existing $C_5,C_6>0$ with
\begin{equation}\label{e196}
\sup_{\epsilon\in D_{\epsilon_0}}|\omega_{p,n}(u,\epsilon)|\le C_5(C_6)^n\frac{n!}{q^{h^2\Delta}}|u|,
\end{equation}
for all $u\in\mathcal{A}_{p,h}$, where $\Delta$ is fixed in (\ref{e2}).
\end{itemize}
\end{prop}
\begin{proof}
The proof of Proposition 5 in~\cite{malek22} can be directly adapted to this parametric framework under the action of the perturbation parameter $\epsilon$. More precisely, the third element in the assumption (\ref{e1}) guarantees that
$$\sup_{\epsilon\in D_{\epsilon_0}}|\epsilon|^{\Delta_{\underline{\ell}}-\ell_0}\le \epsilon_0^{\Delta_{\underline{\ell}}-\ell_0},$$
for every $\underline{\ell}=(\ell_0,\ell_1,\ell_2,\ell_3)\in\mathcal{A}$. The polynomial nature of the coefficients $c_{\underline{\ell}}(z,\epsilon)=\sum_{h\in I_{\underline{\ell}}}c_{\underline{\ell},h}z^h\in \mathcal{O}_b(D_{\epsilon_0})[z]$ allows us to upper estimate
$$\sup_{\epsilon\in D_{\epsilon_0}}|c_{\underline{\ell}}(z,\epsilon)|\le\sum_{h\in I_{\underline{\ell}}}\left(\sup_{\epsilon\in D_{\epsilon_0}}|c_{\underline{\ell},h}(\epsilon)|\right)|z|^h.$$
Taking the above facts into consideration, the proof follows directly from that of Proposition 5 in~\cite{malek22}.
\end{proof}

At this point, one can achieve the construction of actual solutions to the main Cauchy problem (\ref{epral}), (\ref{epralCD}). For that purpose, we introduce the following geometric  definition.

\begin{defin}\label{defi1}
Let $\varsigma\ge 2$ be an integer. Let us consider a finite set of bounded sectors $\underline{\mathcal{E}}=(\mathcal{E}_{p})_{0\le p\le \varsigma-1}$ which conforms a good covering in $\C^{\star}$ (see Definition~\ref{defin-goodcovering}). Associated to $\underline{\mathcal{E}}$ let us choose a bounded sector with vertex at the origin, $\mathcal{T}\subseteq D_{r_{\mathcal{T}}}$ for some $r_{\mathcal{T}}>0$, and a family $\underline{\mathcal{U}}=(\mathcal{U}_p)_{0\le p\le \varsigma-1}$ of unbounded sectors with vertex at the origin, under the following conditions:
\begin{itemize}
\item The roots of $u\mapsto P(ku^k)$ belong to $\C\setminus\left(\bigcup_{p=0}^{\varsigma-1}\mathcal{U}_p\right)$,
\item For every $0\le p\le \varsigma-1$, there exists $\Delta_p>0$ such that for all $t\in\mathcal{T}$ and $\epsilon\in\mathcal{E}_{p}$ a direction $\gamma_p\in\R$ exists (which depends on $t$ and $\epsilon$) such that 
\begin{itemize}
\item[-] $L_{\gamma_p}=[0,\infty)e^{\sqrt{-1}\gamma_p}\subseteq \mathcal{U}_p\cup\{0\}$, and
\item[-] $\cos(k(\gamma_p-\hbox{arg}(\epsilon t)))>\Delta_p$.
\end{itemize}
\end{itemize}
We say that the tuple $\{\mathcal{T},\underline{\mathcal{U}}\}$ is admissible with respect to the good covering $\underline{\mathcal{E}}$.
\end{defin}

\begin{theo}\label{teo1}
Let $\varsigma_1\ge2$. Let $\{\mathcal{T},\underline{\mathcal{U}_1}=(\mathcal{U}_{p,1})_{0\le p\le \varsigma_1-1}\}$ be an admissible set with respect to a given good covering $\underline{\mathcal{E}}=(\mathcal{E}_p)_{0\le p\le \varsigma_1-1}$ (following Definition~\ref{defi1}). Take for granted the assumptions made in Section~\ref{sec21} on the elements involved in the main Cauchy problem (\ref{epral}), (\ref{epralCD}). 

Then, for every $0\le p\le \varsigma_1-1$, there exists a solution $u_{p,1}$ of (\ref{epral}), (\ref{epralCD}) which is holomorphic and bounded on $\mathcal{T}\times D_{R}\times \mathcal{E}_p$. Such solution is of the form
\begin{equation}\label{e224}
u_{p,1}(t,z,\epsilon)=\sum_{n\ge0}u_{p,1,n}(t,\epsilon)\frac{z^n}{n!},\quad (t,z,\epsilon)\in \mathcal{T}\times D_{R}\times \mathcal{E}_p,
\end{equation}
with 
\begin{equation}\label{e225}
u_{p,1,n}(t,\epsilon)=k\int_{L_{\gamma_p}}\omega_{p,1,n}(u,\epsilon)\exp\left(-\left(\frac{u}{\epsilon t}\right)^k\right)\frac{du}{u},
\end{equation}
for every $n\ge0$, and all $(t,\epsilon)\in\mathcal{T}\times\mathcal{E}_p$. Here, $L_{\gamma_p}=[0,\infty)e^{\sqrt{-1}\gamma_p}$, where $\gamma_p$ is determined in Definition~\ref{defi1}, and $\omega_{p,1,n}$ represents the coefficient $\omega_{p,n}$ in (\ref{eA}) of Proposition~\ref{prop5} for the set $\underline{\mathcal{U}}=\underline{\mathcal{U}_1}$.  
\end{theo}
\begin{proof}
Let $0\le p\le \varsigma_1-1$, and consider the function $\omega_p$ constructed in Proposition~\ref{prop5} for the covering $\underline{\mathcal{U}}=\underline{\mathcal{U}_1}$, that we denote $\omega_{p,1}$ with coefficients $\omega_{p,1,n}$. The function $\omega_p$ is a solution of the auxiliary problem (\ref{eaux}), (\ref{eauxCD}). We write $\omega_p$ in the form (\ref{eA}) and recall that its growth with respect to $u$ at infinity, displayed in (\ref{e194}), allows us to apply Laplace transform of order $k$ along direction $\gamma_p$. The construction of the admissible set with respect to the good covering $\underline{\mathcal{E}}$ guarantees that 
$$u_{p,1,n}(t,\epsilon)=\mathcal{L}^{\gamma_p}_{k}(\omega_{p,1,n}(u,\epsilon))(\epsilon t, \epsilon)$$
is holomorphic and bounded in $\mathcal{T}\times\mathcal{E}_p$. Moreover, from the bounds in (\ref{e194}) the expression in (\ref{e224}) determines a bounded holomorphic function in $\mathcal{T}\times D_{R}\times \mathcal{E}_p$. The properties of Laplace transform in Section~\ref{secannex} guarantee that (\ref{e224}) determines an actual solution of (\ref{epral}), (\ref{epralCD}).
\end{proof}

A symmetric construction can be followed by interchanging the role of the variables $\epsilon$ and $t$ in the previous construction to achieve analytic solutions to the main problem in different families of domains. In this respect, Definition~\ref{defi1} reads as follows. 

\begin{defin}\label{defi1bis}
Let $\varsigma\ge 2$ be an integer. Let us consider a finite set of bounded sectors $\underline{\mathcal{T}}=(\mathcal{T}_{p})_{0\le p\le \varsigma-1}$ which conforms a good covering in $\C^{\star}$. Associated to $\underline{\mathcal{T}}$ one chooses a bounded sector with vertex at the origin, $\mathcal{E}\subseteq D_{\epsilon_0}$, and a family $\underline{\mathcal{U}}=(\mathcal{U}_p)_{0\le p\le \varsigma-1}$ of unbounded sectors with vertex at the origin, under the following conditions:
\begin{itemize}
\item The roots of $u\mapsto P(ku^k)$ belong to $\C\setminus\left(\bigcup_{p=0}^{\varsigma-1}\mathcal{U}_p\right)$,
\item For every $0\le p\le \varsigma-1$, there exists $\hat{\Delta}_p>0$ such that for all $\epsilon\in\mathcal{E}$ and $t\in\mathcal{T}_{p}$ a direction $\hat{\gamma}_p\in\R$ exists (which depends on $t$ and $\epsilon$) such that 
\begin{itemize}
\item[-] $L_{\hat{\gamma}_p}=[0,\infty)e^{\sqrt{-1}\hat{\gamma}_p}\subseteq \mathcal{U}_p\cup\{0\}$, and
\item[-] $\cos(k(\hat{\gamma}_p-\hbox{arg}(\epsilon t)))>\hat{\Delta}_p$.
\end{itemize}
\end{itemize}
We say that the tuple $\{\mathcal{E},\underline{\mathcal{U}}\}$ is admissible with respect to the good covering $\underline{\mathcal{T}}$.
\end{defin}

The following symmetric result to Theorem~\ref{teo1} is also valid, following the same arguments as above for its proof.

\begin{theo}\label{teo1bis}
Let $\varsigma_2\ge2$. Let $\{\mathcal{E},\underline{\mathcal{U}_2}=(\mathcal{U}_{p,2})_{0\le p\le \varsigma_2-1}\}$ be an admissible set with respect to a given good covering $\underline{\mathcal{T}}=(\mathcal{T}_p)_{0\le p\le \varsigma_2-1}$ (in the sense of Definition~\ref{defi1bis}). Assume the hypotheses made in Section~\ref{sec21} on the elements involved in the construction of the main Cauchy problem (\ref{epral}), (\ref{epralCD}) hold. Then, for every $0\le p\le \varsigma_2-1$, there exists a solution $u_{p,2}$ of (\ref{epral}), (\ref{epralCD}) which is holomorphic and bounded on $\mathcal{T}_p\times D_{R}\times \mathcal{E}$. Such solution is of the form
\begin{equation}\label{e224bis}
u_{p,2}(t,z,\epsilon)=\sum_{n\ge0}u_{p,2,n}(t,\epsilon)\frac{z^n}{n!},\quad (t,z,\epsilon)\in \mathcal{T}_p\times D_{R}\times \mathcal{E},
\end{equation}
with 
\begin{equation}\label{e225bis}
u_{p,2,n}(t,\epsilon)=k\int_{L_{\hat{\gamma}_p}}\omega_{p,2,n}(u,\epsilon)\exp\left(-\left(\frac{u}{\epsilon t}\right)^k\right)\frac{du}{u},
\end{equation}
for every $n\ge0$, and all $(t,\epsilon)\in\mathcal{T}_p\times\mathcal{E}$. Here, $L_{\hat{\gamma}_p}=[0,\infty)e^{\sqrt{-1}\hat{\gamma}_p}$, is determined in Definition~\ref{defi1bis} and where $\omega_{p,2,n}$ represents the coefficient $\omega_{p,n}$ in (\ref{eA}) of Proposition~\ref{prop5} for the set $\underline{\mathcal{U}}=\underline{\mathcal{U}_2}$. 
\end{theo}

\section{Asymptotic behavior of the analytic solutions, I}\label{secasI}

In this section, we provide two results regarding the asymptotic expansions of the previous analytic solutions determined in Theorem~\ref{teo1} with respect to the perturbation parameter regarding each of the elements in a good covering $\underline{\mathcal{E}}=(\mathcal{E}_{p})_{0\le p\le \varsigma_1-1}$, for some fixed $\varsigma_1\ge2$. The main first result of the present work determines classical Gevrey asymptotic expansions by introducing $L_1-q-$relative-sup-norms on the partial functions $(t,z)\mapsto u_{p,1}(t,z,\epsilon)$, constructed in Theorem~\ref{teo1}, for each fixed $\epsilon\in\mathcal{E}_p$, and all $0\le p\le \varsigma_1-1$, whereas the second result deals with uniform sup-norms on the partial functions $(t,z)$.

\begin{defin}\label{defi243}
let $q>1$ be a real number. We fix a bounded sector with vertex at the origin $\mathcal{T}$. Let us consider the set $\mathcal{O}_b^q(\mathcal{T})[[z]]$ of formal power series $h$ of the form 
$$h(t,z)=\sum_{n\ge0}h_n(t)\frac{z^n}{n!},$$
where $h_n\in\mathcal{O}_b(\mathcal{T}\cap D_{1/q^n})$ for every $n\ge0$.

Let $R_1>0$ be a real number.

We denote by $\mathbb{E}_{L_1;q;R_1}$ the vector space of formal power series $h\in\mathcal{O}_b^{q}(\mathcal{T})[[z]]$ of the form $h(t,z)=\sum_{n\ge0}h_n(t)\frac{z^n}{n!}$, where $h_n\in\mathcal{O}_b(\mathcal{T}\cap D_{1/q^n})$ for every $n\ge0$, such that the $L_1-q-$relative-sup-norm of $h$, defined by
$$\left\|h(t,z)\right\|_{L_1;q;R_1}:=\sum_{n\ge0}\sup_{t\in\mathcal{T}\cap D_{1/q^n}}|h_n(t)|\frac{R_1^n}{n!},$$
is finite.
\end{defin}

\begin{prop}\label{prop0}
The pair $(\mathbb{E}_{L_1;q;R_1},\left\|\cdot\right\|_{L_1;q;R_1})$ is a complex Banach space. The vector space $\mathcal{O}_b(\mathcal{T}\times D_R)$ is contained in $\mathbb{E}_{L_1;q;R_1}$ provided that $R>R_1$.
\end{prop}
\begin{proof}
 First, observe that any formal power series $h(t,z)=\sum_{n\ge0}h_n(t)\frac{z^n}{n!}\in\mathbb{E}_{L_1;q;R_1}$ such that $\left\|h(t,z)\right\|_{L_1;q;R_1}=0$ entails that $h_n\equiv 0$ in $\mathcal{T}\cap D_{1/q^n}$.

Let $(h^p)_{p\ge1}$ be a Cauchy sequence in $\mathbb{E}_{L_1;q;R_1}$. Let us write $h^{p}(t,z)=\sum_{n\ge0}h_n^p(t)\frac{z^n}{n!}$, for every $p\in\N$. It is clear that for every $n\in\N_0$, the sequence $(h^p_n)_{p\ge1}$ is a Cauchy sequence in the Banach space $(\mathcal{O}_b(\mathcal{T}\cap D_{1/q^n}),\left\|\cdot\right\|_{\infty})$, where $\left\|\cdot\right\|_{\infty}$ stands for the supremum norm with 
\begin{equation}\label{e264}
\left\|h^{p}_n-h^{q}_n\right\|_{\infty}<\frac{\epsilon n!}{R_1^n}
\end{equation}
for all $p,q\in\N$ with $p,q\ge p_0$ for some $p_0\in\N$ whenever $\left\|h^p-h^q\right\|_{L_1;q;R_1}<\epsilon$. Therefore, for every $n\in\N_0$ there exists $H_n\in\mathcal{O}_b(\mathcal{T}\cap D_{1/q^n})$ such that $(h^p_n)_{n\ge0}$ converges in sup norm to $H_n$. By considering the formal power series $H(t,z)=\sum_{n\ge0}H_n(t)\frac{z^n}{n!}$ it is direct to check that $H$ is an element of $\mathbb{E}_{L_1;q;R_1}$, with $H$ being the limit of $(h^p)_{p\ge0}$ in such space. 
\end{proof}

From now on, we resume the assumptions and constructions related to the main problem (\ref{epral}), (\ref{epralCD}) in Section~\ref{sec21}. In addition to this, we consider $\varsigma_1\ge2$ and an admissible set $\{\mathcal{T},\underline{\mathcal{U}_1}=(\mathcal{U}_{p,1})_{0\le p\le \varsigma_1-1}\}$ with respect to a given good covering $\underline{\mathcal{E}}=(\mathcal{E}_p)_{0\le p\le \varsigma_1-1}$ (in the sense of Definition~\ref{defi1}). Theorem~\ref{teo1} guarantees the existence of a finite family of solutions $u_{p,1}(t,z,\epsilon)$ to (\ref{epral}), (\ref{epralCD}), for $0\le p\le \varsigma_1-1$, holomorphic and bounded in $\mathcal{T}\times D_R\times \mathcal{E}_p$, of the form (\ref{e224}) with (\ref{e225}), for some $R>0$. Let us fix $0<R_1<R$.

\begin{prop}\label{prop270}
Under the assumptions of Theorem~\ref{teo1}, for every $0\le p\le \varsigma_1-1$ there exist $A_p, B_p>0$ such that
\begin{equation}\label{e272}
\left\|u_{p+1,1}(t,z,\epsilon)-u_{p,1}(t,z,\epsilon)\right\|_{L_1;q;R_1}\le A_p\exp\left(-\frac{B_p}{|\epsilon|^k}\right),
\end{equation} 
valid for every $\epsilon\in\mathcal{E}_p\cap\mathcal{E}_{p+1}$, where we have identified $u_{\varsigma_1,1}:=u_{0,1}$ and $\mathcal{E}_{\varsigma_1}:=\mathcal{E}_0$, provided that $0<R_1<R$ is small enough.
\end{prop}
\begin{proof}
Let $0\le p\le \varsigma_1-1$ and fix $\epsilon\in\mathcal{E}_{p}\cap\mathcal{E}_{p+1}$. In view of (\ref{e224}) and (\ref{e225}), one can write
$$u_{p+1,1}(t,z,\epsilon)-u_{p,1}(t,z,\epsilon)$$
in the form
\begin{equation}\label{e326}
\sum_{n\ge0}\left(k\int_{L_{\gamma_{p+1}}}\omega_{p+1,1,n}(u,\epsilon)\exp\left(-\left(\frac{u}{\epsilon t}\right)^{k}\right)\frac{du}{u}-k\int_{L_{\gamma_{p}}}\omega_{p,1,n}(u,\epsilon)\exp\left(-\left(\frac{u}{\epsilon t}\right)^{k}\right)\frac{du}{u}\right)\frac{z^n}{n!}.
\end{equation}
One can apply b. in Proposition~\ref{prop5} taking $\varsigma_1$ in place of $\varsigma$ and $\underline{\mathcal{U}_1}=(\mathcal{U}_{p,1})_{0\le p\le \varsigma_1-1}$ in place of $\underline{\mathcal{U}}=(\mathcal{U}_{p})_{0\le p\le \varsigma-1}$. Then, for all $n\ge0$, the functions $u\mapsto\omega_{p,1,n}(u,\epsilon)$ and $u\mapsto\omega_{p+1,1,n}(u,\epsilon)$ have a common analytic continuation, say $u\mapsto \omega_n(u,\epsilon)$ on $D_{R_n}$, with $R_n=R_0/q^n$, for some $R_0>0$. This entails that the integration path appearing in the previous difference can be deformed by the application of Cauchy theorem. Hence, $u_{p+1,1}(t,z,\epsilon)-u_{p,1}(t,z,\epsilon)$ equals
\begin{multline}
=\sum_{n\ge0}\left(k\int_{L_{\gamma_{p+1},R_{n+1}}}\omega_{p+1,1,n}(u,\epsilon)\exp\left(-\left(\frac{u}{\epsilon t}\right)^{k}\right)\frac{du}{u}-k\int_{L_{\gamma_{p},R_{n+1}}}\omega_{p,1,n}(u,\epsilon)\exp\left(-\left(\frac{u}{\epsilon t}\right)^{k}\right)\frac{du}{u}\right.\\
\left.+k\int_{C_{\gamma_p,\gamma_{p+1},R_{n+1}}}\omega_{n}(u,\epsilon)\exp\left(-\left(\frac{u}{\epsilon t}\right)^{k}\right)\frac{du}{u}\right)\frac{z^n}{n!},\label{e332}
\end{multline}
where $L_{\gamma_{j},R_{n+1}}=[R_{n+1},\infty)e^{\sqrt{-1}\gamma_{j}}$ for $j=p,p+1$ and $C_{\gamma_p,\gamma_{p+1},R_{n+1}}$ stands for the oriented arc of circle centered at the origin, of radius $R_{n+1}$, from the point $R_{n+1}e^{\sqrt{-1}\gamma_p}$ to $R_{n+1}e^{\sqrt{-1}\gamma_{p+1}}$. 
Let us provide upper bounds for the previous elements.
Let us define
$$I_1(\gamma_{p+1}):=\left|k\int_{L_{\gamma_{p+1},R_{n+1}}}\omega_{p+1,1,n}(u,\epsilon)\exp\left(-\left(\frac{u}{\epsilon t}\right)^{k}\right)\frac{du}{u} \right|.$$
In view of statement a. in Proposition~\ref{prop5} and the fact that $\{\mathcal{T},\underline{\mathcal{U}_1}\}$ is admissible with respect to $\underline{\mathcal{E}}$, one has that $I_1(\gamma_{p+1})$ is upper bounded by
\begin{multline*}
kC_3\frac{1}{R^n}n!\int_{R_{n+1}}^{\infty}\exp\left(k_1\log^2(r+u_0)+\alpha\log(r+u_0)\right)\exp\left(-\left(\frac{r}{|\epsilon t|}\right)^k\cos(k(\gamma_{p+1}-\hbox{arg}(\epsilon t)))\right) dr\\
\le kC_3\frac{1}{R^n}n!\int_{R_{n+1}}^{\infty}\exp\left(k_1\log^2(r+u_0)+\alpha\log(r+u_0)\right)\exp\left(-\frac{1}{2}\left(\frac{r}{|\epsilon t|}\right)^k\Delta_{p+1}\right)\\
\hfill\times\exp\left(-\frac{1}{2}\left(\frac{r}{|\epsilon t|}\right)^k\Delta_{p+1}\right) dr\\
\le kC_3\frac{1}{R^n}n!\exp\left(-\frac{1}{2}\left(\frac{R_{n+1}}{|\epsilon t|}\right)^k\Delta_{p+1}\right)\int_{R_{n+1}}^{\infty}\exp\left(k_1\log^2(r+u_0)+\alpha\log(r+u_0)\right)\\
\hfill\times\exp\left(-\frac{1}{2}\left(\frac{r}{|\epsilon t|}\right)^k\Delta_{p+1}\right) dr\\
\le kC_3J_1(\Delta_{p+1})\frac{1}{R^n}n!\exp\left(-\frac{1}{2}\left(\frac{R_{n+1}}{|\epsilon t|}\right)^k\Delta_{p+1}\right),
\end{multline*}
with 
$$J_1(\Delta_{p+1})=\int_0^{\infty}\exp\left(k_1\log^2(r+u_0)+\alpha\log(r+u_0)\right)\exp\left(-\frac{1}{2}\left(\frac{r}{\epsilon_0r_{\mathcal{T}}}\right)^k\Delta_{p+1}\right) dr<\infty.$$
Taking into account that $R_{n+1}=R_0/q^{n+1}$ we conclude from the previous expression that
\begin{equation}\label{e305}
\sup_{t\in\mathcal{T}\cap D_{1/q^n}}I_1(\gamma_{p+1})\le kC_3J_1(\Delta_{p+1})\frac{1}{R^n}n!\exp\left(-\frac{1}{2}\left(\frac{R_{0}}{q|\epsilon|}\right)^k\Delta_{p+1}\right).
\end{equation}
An analogous reasoning yields
$$I_2(\gamma_{p}):=\left|k\int_{L_{\gamma_{p},R_{n+1}}}\omega_{p,1,n}(u,\epsilon)\exp\left(-\left(\frac{u}{\epsilon t}\right)^{k}\right)\frac{du}{u} \right|$$
is upper estimated
\begin{equation}\label{e306}
\sup_{t\in\mathcal{T}\cap D_{1/q^n}}I_2(\gamma_{p})\le kC_3J_1(\Delta_{p})\frac{1}{R^n}n!\exp\left(-\frac{1}{2}\left(\frac{R_{0}}{q|\epsilon|}\right)^k\Delta_{p}\right).
\end{equation}
Finally, we consider
$$I_3(\gamma_{p},\gamma_{p+1}):=\left|k\int_{C_{\gamma_p,\gamma_{p+1},R_{n+1}}}\omega_{n}(u,\epsilon)\exp\left(-\left(\frac{u}{\epsilon t}\right)^{k}\right)\frac{du}{u} \right|.$$
Regarding statement b. in Proposition~\ref{prop5} and the parametrization $u(\theta)=R_0/q^{n+1}e^{i\theta}$ for $\theta$ between $\gamma_{p}$ and $\gamma_{p+1}$, one arrives at the existence of $\tilde{\Delta}_{p,p+1}>0$ such that
\begin{multline}
I_3(\gamma_{p},\gamma_{p+1})\le k\int_{\gamma_p}^{\gamma_{p+1}}C_1(C_2)^n\frac{n!}{q^{n^2\Delta}}\frac{R_0}{q^{n+1}}\exp\left(-\left(\frac{R_0/q^{n+1}}{|\epsilon t|}\right)^k\cos(k(\theta-\hbox{arg}(\epsilon t)))\right)d\theta\\
\le k|\gamma_{p+1}-\gamma_{p}|C_1(C_2)^n\frac{n!R_0}{q^{n^2\Delta+n+1}}\exp\left(-\left(\frac{R_0/q^{n+1}}{|\epsilon t|}\right)^k\tilde{\Delta}_{p,p+1}\right).\label{e364}
\end{multline}
This entails 
\begin{equation}\label{e307}
\sup_{t\in\mathcal{T}\cap D_{1/q^n}}I_3(\gamma_{p},\gamma_{p+1})\le k|\gamma_{p+1}-\gamma_{p}|C_1(C_2)^n\frac{n!R_0}{q^{n^2\Delta+n+1}}\exp\left(-\left(\frac{R_0}{q|\epsilon|}\right)^k\tilde{\Delta}_{p,p+1}\right).
\end{equation}
Collecting the information in (\ref{e305}), (\ref{e306}) and (\ref{e307}), one arrives at 
\begin{multline*}
\left\|u_{p+1,1}(t,z,\epsilon)-u_{p,1}(t,z,\epsilon)\right\|_{L_1;q;R_1}\\
\le \sum_{n\ge0}\left(2kC_3\max\{J_1(\Delta_{p}),J_1(\Delta_{p+1})\}\frac{1}{R^n}\exp\left(-\frac{1}{2}\left(\frac{R_{0}}{q|\epsilon|}\right)^k\min\{\Delta_{p},\Delta_{p+1}\}\right)\right.\\
\left.+k|\gamma_{p+1}-\gamma_{p}|C_1(C_2)^n\frac{R_0}{q^{n^2\Delta+n+1}}\exp\left(-\left(\frac{R_0}{q|\epsilon|}\right)^k\tilde{\Delta}_{p,p+1}\right)\right)R_1^n.
\end{multline*}
The choice of $R_1<R$ yields (\ref{e272}) with 
$$A_p:=\left(\sum_{n\ge0}\left(\frac{R_1}{R}\right)^n+\sum_{n\ge0}\frac{(C_2R_1)^n}{q^{n^2\Delta+n+1}}\right)\left(2kC_3\max\{J_1(\Delta_{p}),J_1(\Delta_{p+1})\}+k|\gamma_{p+1}-\gamma_{p}|C_1 R_0\right),$$
and 
$$B_p:=\frac{R_0^k}{q^k}\min\left\{\frac{1}{2}\min\{\Delta_{p},\Delta_{p+1}\},\tilde{\Delta}_{p,p+1}\right\}.$$
\end{proof}

Proposition~\ref{prop270} guarantees that Theorem~\ref{rst} (RS) (see Section~\ref{secRS}) can be applied when considering the next Banach space of functions.


\begin{theo}\label{teo2}
Under the assumptions of Theorem~\ref{teo1}, there exists a formal power series $\hat{u}_1(t,z,\epsilon)=\sum_{n\ge0}u_{1,n}(t,z)\frac{\epsilon^n}{n!}\in\mathbb{E}_{L_1;q;R_1}[[\epsilon]]$ 
satisfying that there exist $C,M>0$ such that for all $N\ge0$ one has
\begin{equation}\label{e344}
\left\|u_{p,1}(t,z,\epsilon)-\sum_{n=0}^{N}u_{1,n}(t,z)\frac{\epsilon^n}{n!}\right\|_{L_1;q;R_1}\le C M^{N+1}\Gamma\left(\frac{N+1}{k}\right)|\epsilon|^{N+1},
\end{equation}
for all $0\le p\le \varsigma_1-1$ and all $\epsilon\in\mathcal{E}_p$. In other words, $\epsilon\mapsto \hat{u}_1(t,z,\epsilon)$ is the common Gevrey asymptotic expansion of order $1/k$ of the analytic solution $\epsilon\mapsto u_{p,1}(t,z,\epsilon)$, as a function with values in $\mathbb{E}_{L_1;q;R_1}$ in $\mathcal{E}_p$, for all $0\le p\le \varsigma_1-1$.
\end{theo}
\begin{proof}
For all $0\le p\le \varsigma_1-1$, consider the function $u_{p,1}(t,z,\epsilon)$ constructed in Theorem~\ref{teo1}, and define the map  $G_p(\epsilon)$ by $\epsilon\mapsto u_{p,1}(t,z,\epsilon)$ for $\epsilon\in\mathcal{E}_p$, which is viewed as a function with values in the Banach space $(\mathbb{E}_{L_1;q;R_1},\left\|\cdot\right\|_{L_1;q;R_1})$. Observe that one can apply Theorem~\ref{rst} (RS) in view of Proposition~\ref{prop270}, in order to achieve the existence of $\hat{u}_1(t,z,\epsilon)\in  \mathbb{E}_{L_1;q;R_1}[[\epsilon]]$ with the required properties.

\end{proof}

In the second part of this section we provide a different asymptotic expansion with respect to the perturbation parameter, under the action of a different norm compared to that in the first part.  The main result of this section guarantees $q$-Gevrey asymptotic expansions relating the formal and the analytic solutions by incorporating the classical $L_1$-sup norm as follows. The procedure rests on the approach in~\cite{malek22}.

\begin{defin}\label{defi410}
let $q>1$ be a real number. We fix a bounded sector with vertex at the origin $\mathcal{T}$. Let us consider the set $\mathcal{O}_b(\mathcal{T})[[z]]$ of formal power series $h$ of the form 
$$h(t,z)=\sum_{n\ge0}h_n(t)\frac{z^n}{n!},$$
where $h_n\in\mathcal{O}_b(\mathcal{T})$ for every $n\ge0$.

Let $R_1>0$ be a real number.

We denote by $\mathbb{E}_{L_1;R_1}$ the vector space of formal power series $h\in\mathcal{O}_b(\mathcal{T})[[z]]$ of the form $h(t,z)=\sum_{n\ge0}h_n(t)\frac{z^n}{n!}$, where $h_n\in\mathcal{O}_b(\mathcal{T})$ for every $n\ge0$, such that the $L_1-$sup-norm of $h$, defined by
$$\left\|h(t,z)\right\|_{L_1;R_1}:=\sum_{n\ge0}\sup_{t\in\mathcal{T}}|h_n(t)|\frac{R_1^n}{n!},$$
is finite.
\end{defin}

The proof of the next result is analogous to that of Proposition~\ref{prop0}, so it is omitted.

\begin{prop}\label{prop00}
The pair $(\mathbb{E}_{L_1;R_1},\left\|\cdot\right\|_{L_1;R_1})$ is a complex Banach space. The vector space $\mathcal{O}_b(\mathcal{T}\times D_R)$ is contained in $\mathbb{E}_{L_1;R_1}$ provided that $R>R_1$.
\end{prop}

\noindent\textbf{Remark:}  Observe that $\mathcal{O}_b(\mathcal{T})[[z]]\subseteq \mathcal{O}_b^q(\mathcal{T})[[z]]$ for all $q>1$. Furthermore, $\mathbb{E}_{L_1;R_1}\subseteq\mathbb{E}_{L_1;q;R_1}$ for all $q>1$ and we have the next inequality $\left\|h\right\|_{L_1;q;R_1}\le \left\|h\right\|_{L_1;R_1}$ for all $h\in\mathbb{E}_{L_1;R_1}$. 

\vspace{0.3cm}

In this section, we resume the assumptions and constructions associated to (\ref{epral}), (\ref{epralCD}) of Section~\ref{sec21}. We also fix a good covering $\underline{\mathcal{E}}=(\mathcal{E}_p)_{0\le p\le\varsigma_1-1}$ for some $\varsigma_1\ge2$, and an admissible set (in the sense of Definition~\ref{defi1}) $\{\mathcal{T},\underline{\mathcal{U}_1}=(\mathcal{U}_{p,1})_{0\le p\le \varsigma_1-1}\}$ associated to the previous good covering. Let $u_{p,1}(t,z,\epsilon)\in\mathcal{O}_b(\mathcal{T}\times D_R\times\mathcal{E}_p)$ be the analytic solution to (\ref{epral}), (\ref{epralCD}) for $0\le p\le \varsigma_1-1$, for some fixed $R>0$. We also choose $0<R_1<R$.

\begin{prop}\label{prop422}
Under the assumptions of Theorem~\ref{teo1}, for every $0\le p\le \varsigma_1-1$ and $0<R_1<R$, there exist $\tilde{A}_p,\tilde{B}_p>0$ such that
$$\left\|u_{p+1,1}(t,z,\epsilon)-u_{p,1}(t,z,\epsilon)\right\|_{L_1;R_1}\le\tilde{A}_p(\tilde{B}_p)^N\Gamma\left(\frac{N}{k}\right)q^{N^2/2}|\epsilon|^{N},$$
for every $\epsilon\in\mathcal{E}_p\cap\mathcal{E}_{p+1}$, where we have identified $u_{\varsigma_1,1}:=u_{0,1}$ and $\mathcal{E}_{\varsigma_1}:=\mathcal{E}_0$.
\end{prop}
\begin{proof}
The proof of Theorem 1~\cite{malek22} can be followed point by point, by considering uniform bounds with respect to $t\in\mathcal{T}$, and by assuming that $\epsilon$ plays the role of $t$ variable in Theorem 1 of~\cite{malek22}, together with the estimates in Proposition~\ref{prop5}. 
\end{proof}

From the previous result, one achieves the following second asymptotic relation of Gevrey mixed order $(1/k;(q,1))$. 

\begin{theo}\label{teo3}
Under the assumptions of Theorem~\ref{teo1}, there exists a formal power series $\hat{\tilde{u}}_1(t,z,\epsilon)=\sum_{n\ge0}\tilde{u}_{1,n}(t,z)\frac{\epsilon^n}{n!}\in\mathbb{E}_{L_1;R_1}[[\epsilon]]$ 
satisfying that there exist $C,M>0$ such that for all $N\ge0$ one has
\begin{equation}\label{e344b}
\left\|u_{p,1}(t,z,\epsilon)-\sum_{n=0}^{N}\tilde{u}_{1,n}(t,z)\frac{\epsilon^n}{n!}\right\|_{L_1;R_1}\le C M^{N+1}\Gamma\left(\frac{N+1}{k}\right)q^{\frac{(N+1)^2}{2}}|\epsilon|^{N+1},
\end{equation}
for all $0\le p\le \varsigma_1-1$ and all $\epsilon\in\mathcal{E}_p$. 
\end{theo}
\begin{proof}
One can define the map $G_p(\epsilon)$ for $0\le p\le \varsigma_1-1$ as in the proof of Theorem~\ref{teo2} and apply  Theorem~\ref{trsq} (q-RS) (see Section~\ref{secRS}) in virtue of Proposition~\ref{prop422}.
\end{proof}

\section{Asymptotic behavior of the analytic solutions, II}\label{secasII}

In this section, we describe the asymptotic behavior of the analytic solutions of (\ref{epral}), (\ref{epralCD}) with respect to the time variable, near the origin. We proceed in a similar manner to Section~\ref{secasI}, by fixing an element in a good covering $\underline{\mathcal{T}}$, and providing the asymptotic properties regarding different norms on the partial sums $(z,\epsilon)\mapsto u_{p,2}(t,z,\epsilon)$, for each $t\in\mathcal{T}_p$ and all $0\le p\le \varsigma_2-1$: an $L_1-q-$relative-sup norm, and the uniform sup-norm on the partial functions $(z,\epsilon)$.

In accordance to Definition~\ref{defi243}, one can state the following symmetric definition.

\begin{defin}\label{defi243bis}
let $q>1$ be a real number. We fix a bounded sector with vertex at the origin $\mathcal{E}$. Let us consider the set $\mathcal{O}_b^q(\mathcal{E})[[z]]$ of formal power series $g$ of the form 
$$g(z,\epsilon)=\sum_{n\ge0}g_n(\epsilon)\frac{z^n}{n!},$$
where $g_n\in\mathcal{O}_b(\mathcal{E}\cap D_{1/q^n})$ for every $n\ge0$.

Let $R_1>0$ be a real number.

We denote by $\mathbb{F}_{L_1;q;R_1}$ the vector space of formal power series $g\in\mathcal{O}_b^{q}(\mathcal{E})[[z]]$ of the form $g(z,\epsilon)=\sum_{n\ge0}g_n(\epsilon)\frac{z^n}{n!}$, where $g_n\in\mathcal{O}_b(\mathcal{E}\cap D_{1/q^n})$ for every $n\ge0$, such that the $L_1-q-$relative-sup-norm of $g$, defined by
$$\left\|g(z,\epsilon)\right\|_{L_1;q;R_1}:=\sum_{n\ge0}\sup_{\epsilon\in\mathcal{E}\cap D_{1/q^n}}|g_n(\epsilon)|\frac{R_1^n}{n!},$$
is finite.
\end{defin}

\vspace{0.3cm}

\textbf{Remark:} We adopt the same notation for the norm in Definition~\ref{defi243} for simplicity. The pair $(\mathbb{F}_{L_1;q;R_1},\left\|\cdot\right\|_{L_1;q;R_1})$ is a complex Banach space. Besides, the vector space $\mathcal{O}_b(D_R\times\mathcal{E})$ is contained in $\mathbb{F}_{L_1;q;R_1}$ provided that $R>R_1$.

\vspace{0.3cm}

Returning to the main problem (\ref{epral}), (\ref{epralCD}), we take for granted the assumptions made on it in Section~\ref{sec21}. Let us first fix $\varsigma_2\ge2$ and a family of admissible domains which determines the geometry of the problem. We choose a finite family of sectors $\underline{\mathcal{T}}=(\mathcal{T}_{p})_{0\le p\le \varsigma_2-1}$, which determines a good covering in $\C^{\star}$ (see Definition~\ref{defin-goodcovering}). Associated to such good covering, we consider a tuple $\{\mathcal{E},\underline{\mathcal{U}_2}=(\mathcal{U}_{p,2})_{0\le p\le \varsigma_2-1}\}$ which is admissible with respect to the good covering $\underline{\mathcal{T}}$ (see Definition~\ref{defi1bis}). We assume $\mathcal{E}\subseteq D_{\epsilon_0}$. Theorem~\ref{teo1bis} provides with a family of holomorphic solutions $u_{p,2}(t,z,\epsilon)$, for $0\le p\le \varsigma_2-1$, of the main problem (\ref{epral}), (\ref{epralCD}), which are holomorphic and bounded in $\mathcal{T}_{p}\times D_R\times \mathcal{E}$, for some $R>0$. In addition to this, the form of such solutions is determined by (\ref{e224bis}) with (\ref{e225bis}) for all $n\ge0$. Let $0<R_1<R$.

\begin{prop}\label{prop270bis}
We take for granted the assumptions made in Theorem~\ref{teo1bis}. Then, for every $0\le p\le \varsigma_2-1$ there exist $\tilde{A}_p,\tilde{B}_p>0$ such that
$$\left\|u_{p+1,2}(t,z,\epsilon)-u_{p,2}(t,z,\epsilon)\right\|_{L_1;q;R_1}\le \tilde{A}_p\exp\left(-\frac{\tilde{B}_p}{|t|^k}\right),$$
for every $t\in\mathcal{T}_{p}\cap\mathcal{T}_{p+1}$, where $u_{\varsigma_2,2}$ and $\mathcal{T}_{\varsigma_2}$ stand for $u_{0,2}$ and $\mathcal{T}_{0}$, respectively, provided that $R_1>0$ is small enough.
\end{prop}
\begin{proof}
The proof of this result heavily rests on that of Proposition~\ref{prop270}. We only give details on the steps in which this proof differs from that. 

Let us fix $0\le p\le \varsigma_2-1$ and $t\in\mathcal{T}_p\cap\mathcal{T}_{p+1}$. We write the difference of two consecutive solutions $u_{p+1,2}(t,z,\epsilon)-u_{p,2}(t,z,\epsilon)$ in the form (\ref{e326}), with $L_{\hat{\gamma}_{j}}$ in place of $L_{\gamma_j}$ for $j\in\{p,p+1\}$. Once again, one can apply b. in Proposition~\ref{prop5}, by substituting $\varsigma$ by $\varsigma_2$ and $\underline{\mathcal{U}_2}=(\mathcal{U}_{p,2})_{0\le p\le \varsigma_2-1}$ replacing $\underline{\mathcal{U}}=(\mathcal{U}_{p})_{0\le p\le \varsigma-1}$. For each $n\ge0$, an analogous analytic continuation $u\mapsto \omega_n(u,\epsilon)$ of the functions $u\mapsto\omega_{p,2,n}(u,\epsilon)$ and $u\mapsto\omega_{p+1,2,n}(u,\epsilon)$ for all $\epsilon\in\mathcal{E}$ allows us to rewrite the difference of the consecutive solutions in the form of (\ref{e332}), with
 $L_{\hat{\gamma}_{j},R_{n+1}}$ and $C_{\hat{\gamma}_p,\hat{\gamma}_{p+1},R_{n+1}}$ instead of $L_{\gamma_{j},R_{n+1}}$ and $C_{\gamma_p,\gamma_{p+1},R_{n+1}}$, for $j\in\{p,p+1\}$. From the proof of Proposition~\ref{prop270}, we derive 
$$I_1(\hat{\gamma}_{p+1})\le k C_3J_1(\hat{\Delta}_{p+1})\frac{1}{R^n}n!\exp\left(-\frac{1}{2}\left(\frac{R_{n+1}}{|\epsilon t|}\right)^{k}\hat{\Delta}_{p+1}\right),$$
leading to 
\begin{equation}\label{e305bis}
\sup_{\epsilon\in\mathcal{E}\cap D_{1/q^n}} I_1(\hat{\gamma}_{p+1})\le k C_3J_1(\hat{\Delta}_{p+1})\frac{1}{R^n}n!\exp\left(-\frac{1}{2}\left(\frac{R_{0}}{q |t|}\right)^{k}\hat{\Delta}_{p+1}\right).
\end{equation}
It is straight to achieve the bounds
\begin{equation}\label{e306bis}
\sup_{\epsilon\in\mathcal{E}\cap D_{1/q^n}} I_2(\hat{\gamma}_{p})\le k C_3J_1(\hat{\Delta}_{p})\frac{1}{R^n}n!\exp\left(-\frac{1}{2}\left(\frac{R_{0}}{q |t|}\right)^{k}\hat{\Delta}_{p}\right).
\end{equation}
Finally, taking (\ref{e364}) into account, we arrive at
$$I_3(\hat{\gamma}_p,\hat{\gamma}_{p+1})\le k|\hat{\gamma}_{p+1}-\hat{\gamma}_{p}|C_1(C_2)^n\frac{n!R_0}{q^{n^2\Delta+n+1}}\exp\left(-\left(\frac{R_0/q^{n+1}}{|\epsilon t|}\right)^k\tilde{\hat{\Delta}}_{p,p+1}\right),$$
for some $\tilde{\hat{\Delta}}_{p,p+1}>0$, which leads us to
\begin{equation}\label{e307bis}
\sup_{\epsilon\in\mathcal{E}\cap D_{1/q^n}} I_3(\hat{\gamma}_p,\hat{\gamma}_{p+1})\le k|\hat{\gamma}_{p+1}-\hat{\gamma}_{p}|C_1(C_2)^n\frac{n!R_0}{q^{n^2\Delta+n+1}}\exp\left(-\left(\frac{R_0}{q|t|}\right)^k\tilde{\hat{\Delta}}_{p,p+1}\right).
\end{equation}
In view of (\ref{e305bis}), (\ref{e306bis}) and (\ref{e307bis}) we have

\begin{multline*}
\left\|u_{p+1,2}(t,z,\epsilon)-u_{p,2}(t,z,\epsilon)\right\|_{L_1;q;R_1}\\
\le \sum_{n\ge0}\left(2kC_3\max\{J_1(\hat{\Delta}_{p}),J_1(\hat{\Delta}_{p+1})\}\frac{1}{R^n}\exp\left(-\frac{1}{2}\left(\frac{R_{0}}{q|t|}\right)^k\min\{\hat{\Delta}_p,\hat{\Delta}_{p+1}\}\right)\right.\\
\left.+k|\hat{\gamma}_{p+1}-\hat{\gamma}_{p}|C_1(C_2)^n\frac{R_0}{q^{n^2\Delta+n+1}}\exp\left(-\left(\frac{R_0}{q|t|}\right)^k\tilde{\hat{\Delta}}_{p,p+1}\right)\right)R_1^n.
\end{multline*}
The choice of $R_1<R$ and 
$$\tilde{A}_p:=\left(\sum_{n\ge0}\left(\frac{R_1}{R}\right)^n+\sum_{n\ge0}\frac{(C_2R_1)^n}{q^{n^2\Delta+n+1}}\right)\left(2kC_3\max\{J_1(\hat{\Delta}_{p}),J_1(\hat{\Delta}_{p+1})\}+k|\hat{\gamma}_{p+1}-\hat{\gamma}_{p}|C_1 R_0\right),$$
and 
$$\tilde{B}_p:=\frac{R_0^k}{q^k}\min\left\{\frac{1}{2}\min\{\hat{\Delta}_p,\hat{\Delta}_{p+1}\},\tilde{\hat{\Delta}}_{p,p+1}\right\}$$
allows us to conclude the result.
\end{proof}


The following result is a consequence of the classical Ramis-Sibuya Theorem~\ref{rst} (RS), whose proof can be adapted from that of Theorem~\ref{teo2} to the good covering $\underline{\mathcal{T}}$.

\begin{theo}\label{teo4}
Under the assumptions of Theorem~\ref{teo1bis}, there exists a formal power series $\hat{u}_2(t,z,\epsilon)=\sum_{n\ge0}u_{2,n}(z,\epsilon)\frac{t^n}{n!}\in\mathbb{F}_{L_1;q;R_1}[[t]]$ 
satisfying that there exist $\tilde{C},\tilde{M}>0$ such that for all $N\ge0$ one has
\begin{equation}\label{e344bis}
\left\|u_{p,2}(t,z,\epsilon)-\sum_{n=0}^{N}u_{2,n}(z,\epsilon)\frac{t^n}{n!}\right\|_{L_1;q;R_1}\le \tilde{C} \tilde{M}^{N+1}\Gamma\left(\frac{N+1}{k}\right)|t|^{N+1},
\end{equation}
for all $0\le p\le \varsigma_2-1$ and all $t\in\mathcal{T}_p$. Equivalently, the formal power series $t\mapsto \hat{u}_2(t,z,\epsilon)$ is the common Gevrey asymptotic expansion of order $1/k$ of the analytic solution $t\mapsto u_{p,2}(t,z,\epsilon)$, as a function with values in $\mathbb{F}_{L_1;q;R_1}$ in $\mathcal{T}_p$, for all $0\le p\le \varsigma_2-1$.
\end{theo}

As a final step in the asymptotic study of the solutions to the main problem, we describe the results obtained when considering the classical $L_1$-sup norm, in accordance with~\cite{malek22}. The symmetric version of Definition~\ref{defi410} reads as follows.

\begin{defin}\label{defi243bis2}
let $q>1$ be a real number. We fix a bounded sector with vertex at the origin $\mathcal{E}$. Let us consider the set $\mathcal{O}_b(\mathcal{E})[[z]]$ of formal power series $g$ of the form 
$$g(z,\epsilon)=\sum_{n\ge0}g_n(\epsilon)\frac{z^n}{n!},$$
where $g_n\in\mathcal{O}_b(\mathcal{E})$ for every $n\ge0$.

Let $R_1>0$ be a real number.

We denote by $\mathbb{F}_{L_1;R_1}$ the vector space of formal power series $g\in\mathcal{O}_b(\mathcal{E})[[z]]$ of the form $g(z,\epsilon)=\sum_{n\ge0}g_n(\epsilon)\frac{z^n}{n!}$, where $g_n\in\mathcal{O}_b(\mathcal{E})$ for every $n\ge0$, such that the $L_1$-sup-norm of $g$, defined by
$$\left\|g(z,\epsilon)\right\|_{L_1;R_1}:=\sum_{n\ge0}\sup_{\epsilon\in\mathcal{E}}|g_n(\epsilon)|\frac{R_1^n}{n!},$$
is finite.
\end{defin}


As a result, one can follow the proof of Theorem 1~\cite{malek22}, assuming uniform bounds with respect to $\epsilon\in\mathcal{E}$, to arrive at the following result.

\begin{prop}\label{prop541}
Under the assumptions of Theorem~\ref{teo1bis}, for every $0\le p\le \varsigma_2-1$ and $0<R_1<R$, there exist $\tilde{A}_p,\tilde{B}_p>0$ such that
$$\left\|u_{p+1,2}(t,z,\epsilon)-u_{p,2}(t,z,\epsilon)\right\|_{L_1;R_1}\le\tilde{A}_p(\tilde{B}_p)^N\Gamma\left(\frac{N}{k}\right)q^{N^2/2}|t|^{N},$$
for every $t\in\mathcal{T}_p\cap\mathcal{T}_{p+1}$, with $u_{\varsigma_2,2}:=u_{0,2}$ and $\mathcal{T}_{\varsigma_2}:=\mathcal{T}_0$.
\end{prop}

Finally, a direct application of  Theorem~\ref{trsq} (q-RS) to Proposition~\ref{prop541} leads us to the last asymptotic result on the behavior of the analytic solutions of the main problem.

 \begin{theo}\label{teo5}
Under the assumptions of Theorem~\ref{teo1bis}, there exists a formal power series $\hat{\tilde{u}}_2(t,z,\epsilon)=\sum_{n\ge0}\tilde{u}_{2,n}(z,\epsilon)\frac{t^n}{n!}\in\mathbb{F}_{L_1;R_1}[[t]]$ 
satisfying that there exist $\tilde{C},\tilde{M}>0$ such that for all $N\ge0$ one has
\begin{equation}\label{e344c}
\left\|u_{p,2}(t,z,\epsilon)-\sum_{n=0}^{N}\tilde{u}_{2,n}(z,\epsilon)\frac{t^n}{n!}\right\|_{L_1;R_1}\le \tilde{C} \tilde{M}^{N+1}\Gamma\left(\frac{N+1}{k}\right)q^{\frac{(N+1)^2}{2}}|t|^{N+1},
\end{equation}
for all $0\le p\le \varsigma_2-1$ and all $t\in\mathcal{T}_p$. 
\end{theo}

\section{Annex I: Laplace transform}\label{secannex}

In this annex, we briefly describe the definition of Laplace transform and its main properties considered in the present work. We omit its proofs which can be found in previous research by the authors such as~\cite{lama17}.

In the whole section, $(\mathbb{E},\left\|\cdot\right\|_{\mathbb{E}})$ stands for a complex Banach space.

\begin{defin}\label{defi127}
Let $k\in\N$ and choose an unbounded sector $S_{d,\delta}:=\{z\in\C^{\star}:|\hbox{arg}(u)-d|<\delta\}$, for some $d\in\R$ and $\delta>0$. Let $f\in\mathcal{O}(S_{d,\delta},\mathbb{E})$, continuous in $S_{d,\delta}\cup\{0\}$, such that two constants $C,K>0$ exist with
\begin{equation}\label{e129}
\left\|f(u)\right\|_{\mathbb{E}}\le C|u|\exp\left(K|u|^{k}\right),\quad u\in S_{d,\delta}.
\end{equation}
The Laplace transform of order $k$ of $f$ along direction $d$ is 
$$\mathcal{L}_{k}^{d}(f(u))(t)=k\int_{L_{\gamma}}f(u)\exp\left(-\left(\frac{u}{t}\right)^{k}\right)\frac{du}{u},$$
where $L_{\gamma}=[0,\infty)e^{i\gamma}\subseteq S_{d,\delta}\cup\{0\}$ and the value of $\gamma\in\R$ depends on $t$, with $\cos(k(\gamma-\hbox{arg}(t)))\ge \delta_{1}$ for some $\delta_1>0$.

$\mathcal{L}_k^d(f(u))(t)$ defines a holomorphic and bounded function in $\mathcal{O}(S_{d,\theta,R^{1/k}},\mathbb{E})$, with $S_{d,\theta,R^{1/k}}=D(0,R^{1/k})\cap S_{d,\theta/2}$, where $\frac{\pi}{k}<\theta<\frac{\pi}{k}+2\delta$ and $0<R<\delta_1/K$.
\end{defin}

\begin{prop}\label{propLaplace1}
Let $k\in\N$ and $d\in\R$. Let $f\in\mathcal{O}(\C,\mathbb{E})$ such that there exist $C,K>0$ with 
$$\left\|f(u)\right\|_{\mathbb{E}}\le C|u|\exp\left(K|u|^{k}\right),\quad u\in\C.$$ 
Then, $\mathcal{L}_{k}^{d}(f(u)):D\to\mathbb{E}$ defines a holomorphic function on some neighborhood of the origin, $D$, and it holds that 
$$\left(\mathcal{L}_k^d(f(u))\right)^{(p)}(0)=\Gamma\left(\frac{n}{k}\right)f^{(p)}(0),$$
for every $p\in\N_0$.
\end{prop}

\begin{corol}\label{coroLaplace1}
In the situation of the previous result, one has that Laplace transform of order $k$ of the polynomial $p(u)=\sum_{h=0}^{m}a_hu^h\in\mathbb{E}[u]$ is the polynomial $\mathcal{L}_{k}^{d}(p(u))(t)=\sum_{h=0}^{m}\Gamma\left(\frac{h}{k}\right)a_ht^h\in\mathbb{E}[t]$.
\end{corol}

\begin{prop}\label{propLaplace2}
Under the hypotheses of Definition~\ref{defi127}, the function $u\mapsto ku^kf(u)$  and, for all $m\in\N$, the convolution product $u\mapsto u^m\star_kf(u)$ given by
$$u^m\star_kf(u)=\frac{u^k}{\Gamma\left(\frac{m}{k}\right)}\int_0^{u^{k}}(u^k-s)^{\frac{m}{k}-1}f(s^{\frac{1}{k}})\frac{ds}{s},$$
admit Laplace transform of order $k$ along direction $d$ with
$$\mathcal{L}_k^d(ku^kf(u))(t)=t^{k+1}\partial_t\left(\mathcal{L}_k^d(f(u))(t)\right)$$
and
$$\mathcal{L}_k^d\left(u\mapsto (u^m\star_k f(u))\right)(t)=t^m\mathcal{L}_k^d(f(u))(t),$$
for all $t\in S_{d,\theta,R^{1/k}}$ and $0<R<\delta_1/K$.
\end{prop}

\begin{prop}\label{propLaplace3}
Let $q>1$ and $\delta\ge1$. Under the hypotheses of Definition~\ref{defi127}, the function $u\mapsto f(q^{\delta}u)$ admits Laplace transform of order $k$ along direction $d$ and one has that
$$\mathcal{L}_k^d(f(q^{\delta}u))(t)=\mathcal{L}_k^d(f(u))(q^{\delta}t),$$
for every $t\in S_{d,\theta,R^{1/k}}$ and $0<R<\delta_1/(Kq^{k\delta})$.
\end{prop}

\section{Annex II: Ramis-Sibuya type Theorems}\label{secRS}

In this section, we recall different versions of Ramis-Sibuya theorem (RS theorem), whose classical statement can be found in~\cite{hssi}, Lemma XI-2-6. Such type of results have been previously successfully applied in the asymptotic study of analytic solutions to many functional problems: a classical RS theorem in the asymptotic study of the solutions to partial differential equations in the complex domain in~\cite{chenlastramalek,lama19} and also a close version to the previous one adapted to multi-index sectors~\cite{family2}, another RS-type theorem adapted to the more general framework of strongly regular sequences in~\cite{lamasa2}, also a $q$-Gevrey version of RS theorem in the framework of $q$-difference-differential equations in~\cite{lama17}, a mixed (Gevrey and $q$-Gevrey) version of RS theorem in the study of $q$-difference-differential problems in~\cite{malek22}.

Let $\varsigma\ge2$ be an integer. In the whole section, $(\mathbb{E},\left\|\cdot\right\|_{\mathbb{E}})$ stands for a complex Banach space. 

\begin{defin}\label{defin-goodcovering}
A finite family $\underline{\mathcal{E}}=(\mathcal{E}_{p})_{0\le p\le \varsigma-1}$ consisting of bounded sectors with vertex at the origin is said to be a good covering in $\C^{\star}$ whenever the following conditions are satisfied:
\begin{itemize}
\item $\mathcal{E}_j\cap\mathcal{E}_{j+1}\neq\emptyset$ for every $0\le j\le \varsigma-1$ (by convention, we define $\mathcal{E}_{\varsigma}:=\mathcal{E}_{0}$).
\item For every three indices $p_1,p_2,p_3\in\{0,\ldots,\varsigma-1\}$ such that $p_i\neq p_j$ for $i,j\in\{1,2,3\}$ with $i\neq j$, then one has $\mathcal{E}_{p_1}\cap \mathcal{E}_{p_2}\cap \mathcal{E}_{p_3}=\emptyset$.
\item There exists a neighborhood of $0\in\C$, say $D$, such that $\bigcup_{p=0}^{\varsigma-1}\mathcal{E}_{p}=D\setminus\{0\}$.
\end{itemize}
\end{defin}

We recall the classical Ramis-Sibuya theorem.

\begin{theo}[RS]\label{rst}
Let $(\mathcal{E}_p)_{0\le p\le \varsigma-1}$ be a good covering in $\C^{\star}$. Let $G_p:\mathcal{E}_p\mapsto\mathbb{E}$ be a holomorphic map for all $0\le p\le \varsigma-1$ such that the following conditions hold:
\begin{itemize}
\item $G_p\in\mathcal{O}_b(\mathcal{E}_p,\mathbb{E})$ for $0\le p\le \varsigma-1$.
\item Given $0\le p\le \varsigma-1$, the cocycle $\Theta_p(\epsilon):=G_{p+1}(\epsilon)-G_p(\epsilon)$ is exponentially flat of order $k$ on $Z_p:=\mathcal{E}_p\cap\mathcal{E}_{p+1}$ (we write $G_{\varsigma}:=G_{0}$), i.e. there exist $A_p,B_p>0$ such that
$$\left\|\Theta(\epsilon)\right\|_{\mathbb{E}}\le A_p\exp\left(-\frac{B_p}{|\epsilon|^k}\right),\qquad \epsilon\in\mathcal{E}_p\cap\mathcal{E}_{p+1}.$$
\end{itemize}
Then, there exists a formal power series, common for every $0\le p\le \varsigma-1$, $\hat{G}(\epsilon)=\sum_{n\ge0}G_n\epsilon^n\in\mathbb{E}[[\epsilon]]$, such that $G_p$ admits $\hat{G}$ as its Gevrey asymptotic expansion of order $1/k$ on $\mathcal{E}_p$, for all $0\le p\le \varsigma-1$, meaning that there exist $C,M>0$ with
$$\left\|G_p(\epsilon)-\sum_{n=0}^{N}G_n\epsilon^n\right\|_{\mathbb{E}}\le CM^{N+1}\Gamma\left(\frac{N+1}{k}\right)|\epsilon|^{N+1},$$
for all $N\ge0$ and $\epsilon\in\mathcal{E}_p$.
\end{theo}

A mixed Gevrey and $q-$Gevrey version of Ramis-Sibuya theorem is also available.

\begin{theo}[$q-$RS]\label{trsq}
Let $(\mathcal{E}_p)_{0\le p\le \varsigma-1}$ be a good covering in $\C^{\star}$. Let $G_p:\mathcal{E}_p\mapsto\mathbb{E}$ be a holomorphic map for all $0\le p\le \varsigma-1$ such that the following conditions hold:
\begin{itemize}
\item $G_p\in\mathcal{O}_b(\mathcal{E}_p,\mathbb{E})$ for $0\le p\le \varsigma-1$.
\item Given $0\le p\le \varsigma-1$, the cocycle $\Theta_p(\epsilon):=G_{p+1}(\epsilon)-G_p(\epsilon)$ satisfies that there exist $A_p,B_p>0$ with
$$\left\|\Theta(\epsilon)\right\|_{\mathbb{E}}\le A_p (B_p)^N\Gamma\left(\frac{N}{k}\right)q^{\frac{N^2}{2}}|\epsilon|^{N},\qquad \epsilon\in\mathcal{E}_p\cap\mathcal{E}_{p+1},$$
valid for all $N\ge1$.
\end{itemize}
Then, there exists a formal power series, common for every $0\le p\le \varsigma-1$, $\hat{G}(\epsilon)=\sum_{n\ge0}G_n\epsilon^n\in\mathbb{E}[[\epsilon]]$, such that $G_p$ admits $\hat{G}$ as its Gevrey asymptotic expansion of mixed order $(1/k;(q,1))$ on $\mathcal{E}_p$, for all $0\le p\le \varsigma-1$, meaning that there exist $C,M>0$ with
$$\left\|G_p(\epsilon)-\sum_{n=0}^{N}G_n\epsilon^n\right\|_{\mathbb{E}}\le CM^{N+1}\Gamma\left(\frac{N+1}{k}\right)q^{\frac{(N+1)^2}{2}}|\epsilon|^{N+1},$$
for all $N\ge0$ and $\epsilon\in\mathcal{E}_p$.
\end{theo} 

\vspace{0.4cm}

\textbf{Aknowledgements:} Both authors are partially supported by the project PID2019-105621GB-I00 of Ministerio de Ciencia e Innovaci\'on, Spain and by Direcci\'on General de Investigaci\'on e Innovaci\'on. The first author is partially supported by Direcci\'on General de Investigaci\'on e Innovaci\'on, Consejer\'ia de Educaci\'on e Investigaci\'on of Comunidad de Madrid, Universidad de Alcal\'a under grant CM/JIN/2021-014, and by Ministerio de Ciencia e Innovaci\'on-Agencia Estatal de Investigaci\'on MCIN/AEI/10.13039/501100011033 and the European Union ``NextGenerationEU''/ PRTR, under grant TED2021-129813A-I00.

This work was completed during the authors' stay in Tokyo participating in workshops at Josai University, Shibaura Institute of Technology and Chiba University. The authors thank professor Haraoka for his kind invitation.


\begin{thebibliography}{99}
\bibitem{carrillolastra} S. A. Carrillo, A. Lastra, $q-$Nagumo norms and formal solutions to singularly perturbed $q-$difference equations,	arXiv:2307.15096 [math.GM] 
\bibitem{chenlastramalek} G. Chen, A. Lastra, S. Malek, \emph{Parametric Gevrey asymptotics in two complex time variables through truncated Laplace transforms.} Adv. Differ. Equ. 2020, 307 (2020). 
\bibitem{drelasmal} T. Dreyfus, A. Lastra, S. Malek, \emph{Multiple-scale analysis for some linear partial $q-$difference and differential equations with holomorphic coefficients}, Advances in Difference Equations, 2019:326, 2019.
\bibitem{hssi} P. Hsieh, Y. Sibuya, \emph{Basic theory of ordinary differential equations}. Universitext. Springer-Verlag, New York, 1999.
\bibitem{lama19} A. Lastra, S. Malek, \emph{On singularly perturbed linear initial value problems with mixed irregular and Fuchsian time singularities,} J. Geom. Anal. 30 (2020), 3872--3922. 
\bibitem{family2} A. Lastra, S. Malek, \emph{On parametric Gevrey asymptotics for some initial value problems in two asymmetric complex time variables,} Results Math. 73 (2018), no. 4, Art. 155, 46 pp.
\bibitem{lama17} A. Lastra, S. Malek, \emph{On multiscale Gevrey and $q-$Gevrey asymptotics for some linear $q-$difference-differential initial value Cauchy problems}. J. Differ. Equ. Appl. 23(8) (2017), 1397--1457. 
\bibitem{lamasa0} A. Lastra, S. Malek, J. Sanz, \emph{On $q-$asymptotics for linear $q-$difference-differential equations with Fuchsian and irregular singularities}. J. Differ. Equations 252 (2012), no. 10, 5185--5216.
\bibitem{lamasa2} A. Lastra, S. Malek, J. Sanz, \emph{Strongly regular multi-level solutions of singularly perturbed linear partial differential equations,} Results Math. 70 (2016), no. 3--4, 581--614.
\bibitem{ma17} S. Malek, \emph{Parametric Gevrey asymptotics for a $q-$analog of some linear initial value problem}, Funkc. Ekvacioj 60 no. 1 (2017), 21--63.
\bibitem{malek22} S. Malek, \emph{Asymptotics and confluence for some linear $q$-difference-differential Cauchy problem}, J. Geom. Anal.  32, No. 3, Paper No. 93, 63 p. (2022). 
\bibitem{prrasp} D. W. Pravica, N. Randriampiry, M. J. Spurr, \emph{On $q-$advanced spherical Bessel functions of the first kind and perturbations of the Haar wavelet.} Appl. Comput. Harmon. Anal. 44, No. 2 (2018), 350-413. 
\bibitem{prrasp1} D. W. Pravica, N. Randriampiry, M. J. Spurr, \emph{Solutions of a class of multiplicatively advanced differential equations.} C. R., Math., Acad. Sci. Paris 356, No. 7 (2018), 776--817. 
\bibitem{tahara} H. Tahara, \emph{$q-$analogues of Laplace and Borel transforms by means of $q-$exponentials}, Ann. Inst. Fourier, 67 (2017), no. 5, 1865--1903.
\bibitem{taya2} H. Tahara, H. Yamazawa, \emph{$q-$analogue of summability of formal solutions of some linear $q-$difference-differential equations}, Opuscula Math. 35 (2015), no. 5, 713--738.
\bibitem{tuyuan} H. Tu, W. Yuan, \emph{Growth of solutions to two systems of $q-$difference differential equations}, Adv. Differ. Equ. 2020, 112 (2020). 
\bibitem{hongyanlianzhong} H Y. Xu, L Z. Yang, H. Wang, \emph{Growth of the solutions of some q-difference differential equations}, Adv. Differ. Equ. 2015, 172 (2015).
\bibitem{ya1} H. Yamazawa, \emph{Holomorphic and singular solutions of $q-$difference-differential equations of Briot-Bouquet type}, Funkcial. Ekvac. 59 (2016), no. 2, 185--197.
\end{thebibliography}
\end{document}